\documentclass[journal]{IEEEtran}

\ifCLASSINFOpdf
\else
\fi

\usepackage{hhline}
\usepackage{placeins}
\usepackage{amsmath}
\usepackage{multicol}
\usepackage{algorithm}
\usepackage[noend]{algorithmic}
\usepackage[dvipsnames]{xcolor}
\usepackage{color}
\usepackage{tikz,pgfplots}
\usepackage{multirow}
\usepackage{url}
\usetikzlibrary{matrix}
\usetikzlibrary{patterns}
\usepgfplotslibrary{groupplots}
\pgfplotsset{compat=newest}
\usetikzlibrary{calc,positioning,shapes,fit,arrows,shadows}
\tikzstyle{line} = [ draw, -latex']
\usetikzlibrary{shapes,arrows}
\usepgfplotslibrary{units}
\DeclareMathOperator{\atantwo}{atan2}
\usepackage{nomencl}
\makenomenclature
\usepackage[version=4]{mhchem}
\usepackage{subcaption}
\usepackage{etoolbox}
\renewcommand\nomgroup[1]{%
  \item[\bfseries
  \ifstrequal{#1}{A}{Sets and indices}{%
  \ifstrequal{#1}{B}{Parameters}{%
  \ifstrequal{#1}{C}{Decision Variables}{}}}%
]}

\newcommand{\cT}{\mathcal{T}}
\newcommand{\cB}{\mathcal{B}}
\newcommand{\cL}{\mathcal{L}}
\newcommand{\cG}{\mathcal{G}}

\usepackage{color}
\usepackage{float}

\begin{document}

\title{Towards Low-carbon Power Networks: Optimal Integration of Renewable Energy Sources and Hydrogen Storage}

\author{Sezen Ece Kayac{\i}k, Albert H. Schrotenboer, Evrim Ursavas, Iris F. A. Vis
\thanks{
S. E. Kayac{\i}k, E. Ursavas, I. F. A. Vis are with the Department of Operations, Faculty of Economics and Business, University of Groningen, Groningen, Netherlands (e-mails: s.e.kayacik, e.ursavas,  i.f.a.vis@rug.nl).

A. H. Schrotenboer is with the School of Industrial Engineering, Eindhoven University of Technology, Eindhoven, Netherlands (e-mail: a.h.schrotenboer@tue.nl).
}
}
\maketitle
\begin{abstract}
This paper proposes a new optimization model and solution method for determining optimal locations and sizing of renewable energy sources and hydrogen storage in a power network. We obtain these strategic decisions based on the multi-period alternating current optimal power (AC OPF) flow problem that considers the uncertainty of renewable output, electricity demand, and electricity prices. We develop a second-order cone programming approach within a Benders decomposition framework to provide globally optimal solutions. To the best of our knowledge, our paper is the first to propose a systematic optimization framework based on AC OPF that jointly analyzes power network, renewable, and hydrogen storage interactions in order to provide optimal locations and sizing decisions of renewables and hydrogen storage. In a test case, we show that the joint integration of renewable sources and hydrogen storage and consideration of the AC OPF model significantly reduces the operational cost of the power network. In turn, our findings can provide quantitative insights to decision-makers on how to integrate renewable sources and hydrogen storage under different settings of the hydrogen selling price, renewable curtailment costs, emission tax prices, and conversion efficiency.
\end{abstract}


\IEEEpeerreviewmaketitle

\mbox{}

\nomenclature[A]{\(\cB\)}{Set of buses, indexed by $i$}
\nomenclature[A]{\(\cL\)}{Set of lines, indexed by $(i,j)$}
\nomenclature[A]{\(\cG\)}{Set of conventional generators, indexed by $g$ }
\nomenclature[A]{\(\cT\)}{Set of time periods, indexed by $t$}
\nomenclature[A]{\(\Omega\)}{Uncertainty set, indexed by $\omega$}
\nomenclature[A]{\(\delta(i)\)}{Set of neighbors for bus $i$}
\nomenclature[B]{\(s^{\text{max}}\)}{Storage capacity in terms of hours}
\nomenclature[B]{\(f^{\text{rate}}\)}{Ratio for fuel cell}
\nomenclature[B]{\(\underline{r_{i}^{\textsc{r}}}, \overline{r_{i}^{\textsc{r}}} \)}{Min and max allowable power ratings of renewable at bus $i$ }
\nomenclature[B]{\(\underline{h_{i}^{\textsc{r}}}, \overline{h_{i}^{\textsc{r}}} \)}{Min and max allowable power ratings of storage at bus $i$ }
\nomenclature[B]{\(R_i^{\text{up}}\)}{Ramp up limit for generator at bus $i$}
\nomenclature[B]{\(R_i^{\text{down}}\)}{Ramp down limit for generator at bus $i$}
\nomenclature[B]{\( C^C \)}{Renewable curtailment cost }
\nomenclature[B]{\( C^R\)}{Cost per MW of renewable installed }
\nomenclature[B]{\( C^H\)}{Cost per MW of hydrogen storage unit }
\nomenclature[B]{\( C^E \)}{Emission tax price }
\nomenclature[B]{\( C^S \)}{Hydrogen selling price }
\nomenclature[B]{\( C^U \)}{Cost of unsupplied loads}
\nomenclature[B]{\(p_{i}^d\)}{Active power load at bus $i$}
\nomenclature[B]{\(q_{i}^d\)}{Reactive power load at bus $i$}
\nomenclature[B]{\(p_{it}^d(\omega)\)}{Active power load at bus $i$, time $t$, scenario $\omega$}
\nomenclature[B]{\(q_{it}^d(\omega)\)}{Reactive power load at at bus $i$, time $t$, scenario $\omega$}
\nomenclature[B]{\(\underline V_i\)}{Lower bound on the voltage magnitude at bus $i$}
\nomenclature[B]{\(\overline V_i\)}{Upper bound on the voltage magnitude at bus $i$}
\nomenclature[B]{\(g_{ii},b_{ii}\)}{Shunt susceptance at bus $i$}
\nomenclature[B]{\(\overline{(p2g)}_{i}\)}{Maximum allowable power-to-gas conversion at bus $i$}
\nomenclature[B]{\(\overline{(g2p)}_{i}\)}{Maximum allowable gas-to-power conversion at bus $i$}
\nomenclature[B]{\(r_{it}(\omega)\)}{Renewable power factor at bus $i$, time $t$, scenario $\omega$}
\nomenclature[B]{\(\overline B\)}{Total investment budget for renewables and storage}
\nomenclature[B]{\(\eta_{g} \)}{Power-to-gas efficiency}
\nomenclature[B]{\(\eta_{p}\)}{Gas-to-power efficiency}
\nomenclature[B]{\(\underline{p}_i,\overline{p}_i\)}{Lower and upper limits of  active output of generator located at bus $i$}
\nomenclature[B]{\(\underline{q}_i,\overline{q}_i\)}{Lower and upper limits of  reactive output of generator located at bus $i$}
\nomenclature[B]{\(G_{ij}\)}{Conductance for line $(i,j)$}
\nomenclature[B]{\(B_{ij}\)}{Susceptance for line $(i,j)$}
\nomenclature[B]{\(\overline S_{ij}\)}{Maximum allowable flow on line $(i,j)$}
\nomenclature[B]{\(\overline\theta_{ij}\)}{Phase angle bound for line $(i,j)$}
\nomenclature[B]{\(\rho_{\omega}\)}{Probability of scenario $\omega$}
\nomenclature[C]{\(\text{\textbar}V_{it}(\omega) \text{\textbar}\)}{Voltage magnitude at bus $i$, time $t$, scenario $\omega$ }
\nomenclature[C]{\(l_{it}(\omega) \)}{Renewable power curtailment at bus $i$, time $t$, scenario $\omega$ }
\nomenclature[C]{\(u_{it}(\omega) \)}{Unsupplied load at bus $i$, time $t$, scenario $\omega$ }
\nomenclature[C]{\(\theta_{it}(\omega)\)}{Phase angle at bus $i$, time $t$, scenario $\omega$}
\nomenclature[C]{\( r_i^{\textsc{b}}\)}{If a renewable energy source is constructed at bus $i$ $r_i^{\textsc{b}}=1$, and otherwise $r_i^{\textsc{b}}=0$ }
\nomenclature[C]{\( h_i^{\textsc{b}}\)}{If a hydrogen storage is constructed at bus $i$ $h_i^{\textsc{b}}=1$, and otherwise $h_i^{\textsc{b}}=0$ }
\nomenclature[C]{\( r_i^{\textsc{r}}\)}{Power rating of renewable at bus $i$, time $t$, scenario $\omega$ }
\nomenclature[C]{\( h_i^{\textsc{r}}\)}{Power rating of storage at bus $i$, time $t$, scenario $\omega$ }
\nomenclature[C]{\((p2g)_{it}(\omega)\)}{Power-to-gas conversion at bus $i$, time $t$, scenario $\omega$ }
\nomenclature[C]{\((g2p)_{it}(\omega)\)}{Gas-to-power conversion at bus $i$, time $t$, scenario $\omega$ }
\nomenclature[C]{\(s_{it}(\omega)\)}{Energy state-of-charge of storage at bus $i$, time $t$, scenario $\omega$  }
\nomenclature[C]{\(h_{it}(\omega)\)}{ Amount of hydrogen sold at bus $i$, time $t$, scenario $\omega$}
\nomenclature[C]{\(p_{it}^g(\omega)\)}{ Active power output at bus $i$, time $t$, scenario $\omega$ }
\nomenclature[C]{\(q_{it}^g(\omega)\)}{Reactive power output at bus $i$, time $t$, scenario $\omega$ }
\nomenclature[C]{\(p_{ijt}(\omega)\)}{Active power flow at line $(i,j)$, time $t$, scenario $\omega$}
\nomenclature[C]{\(q_{ijt}(\omega)\)}{Reactive power flow at line $(i,j)$, time $t$, scenario $\omega$}

\printnomenclature

\section{Introduction}

\IEEEPARstart{T}{o} achieve net-zero emission targets by 2050 \cite{IEA2020}, governments strongly encourage the deployment of renewable energy production to reduce the emissions caused by electricity and heat generation, which currently accounts for 46\% of the increase in global emissions \cite{IEA2022}. However, the increased penetration of renewable energy into power networks disrupts electricity supply-demand matching due to the intermittency and uncertainty of renewable energy output. The use of green hydrogen, i.e., hydrogen generated from renewable sources, is a high-potential solution to this problem. It can be used to store renewable energy to mitigate supply-demand imbalances of electricity. It can also be sold outside the network to satisfy green hydrogen demand from various sectors, including industry and mobility, providing new economic opportunities \cite{li2021value,schrotenboer2022green}. This paper studies how renewables and hydrogen storage can be integrated into existing power networks efficiently. To the best of our knowledge, this is the first study to design an integrated system of a power network, renewables, and hydrogen storage by providing optimal location and sizing decisions of renewables and hydrogen storage.

A power network operator is responsible for ensuring the network's reliability and cost-efficiency at the operational level \cite{tabares2015multistage}. The network's structure regarding the location and sizing of renewables and hydrogen storage significantly affects operational planning. From a technical perspective, improper placement of renewables and storage causes challenges, including high power losses, voltage instability, and power quality and protection degradation \cite{ehsan2018optimal}. From an economic perspective, the high-capital costs of renewables and storage should be worth the resulting daily operational gains. Therefore, it is crucial to determine the strategic location and sizing decisions considering daily network operations to provide a reliable power network and fully exploit the economic and environmental benefits of renewables and hydrogen storage.

We provide a new model together with a solution approach for integrating renewables and storage into power networks while explicitly considering the operational level challenges. In this regard, the literature has provided valuable contributions, but only to isolated parts of this joint optimization problem. We review this literature in four steps. First, we discuss recent studies on integrating renewables, and second, on integrating general storage types. Third, we provide an overview of the recent works on operational level planning, i.e., optimal power flow (OPF). Last, we outline the new characteristics introduced by considering green hydrogen as a storage type in our setting.

The first group of papers studies only renewable integration into power networks; see \cite{ehsan2018optimal} for an overview. Recently, the location of renewables has been studied while ignoring sizing decisions in radial distribution networks \cite{8027196}, for which the authors provide a heuristic approach considering uncertainty in renewable output and network demand. In \cite{8886404}, the same uncertainties are tackled, but only the sizing of renewables is considered. However, both studies neglect joint location and sizing decisions, which may result in suboptimal decisions. In \cite{yazdavar2020optimal}, joint location and sizing decisions are studied considering renewable intermittency. While the aforementioned studies draw conclusions about the integration of renewables into power networks, the need for more accurate and computationally efficient solution methods is emphasized in \cite{ehsan2018optimal}.

The second group of papers studies the integration of energy storage systems into power networks with given locations and sizing of renewables \cite{das2018overview}. In \cite{fernandez2016optimal}, a direct current (DC) OPF model is proposed considering renewable uncertainty to determine locations and sizing of storage systems in a transmission network, aiming to minimize the total operating cost and the investment cost of storage systems. They show that the operational level parameters, such as curtailment cost of renewables, affect the location and sizing decisions. A similar DC OPF optimization study proposed in \cite{7317588} shows that increasing the capital investment in storage systems can reduce the daily operating cost of the power network. The authors of \cite{zheng2019hierarchical} and \cite{9431711} propose hierarchical planning models considering alternating current (AC) power equations for radial distribution networks. These methods, however, cannot be directly applied to meshed transmission networks since power flows frequently change direction throughout the day or as a function of the production from generators \cite{fernandez2016optimal}.  

Next to the need for further developments in storage studies, very limited literature is available on the joint integration of renewables and storage. The authors of \cite{8476295} and \cite{garcia2019optimal} are among the initial attempts for joint optimization. The authors of \cite{ehsan2018coordinated} study locations and sizing of both renewables and energy storage systems. They use the heuristic moment matching method to represent renewable output and network load uncertainties. However, their study does not include interactions with conventional generators and uses a local optimization method. A literature review by \cite{ehsan2018optimal} further emphasizes the need for joint studies since the combined planning of renewables and energy storage systems can increase the reliability and power quality of power networks.  

For efficient location and sizing decisions, the underlying operational level problem needs to be analyzed carefully. Mainly OPF models are used because they can analyze the impact of location and sizing decisions on daily network operations. Most studies consider a 24-hour horizon due to hourly fluctuations in demand and supply; however, considering such a long horizon poses a computational burden. Albeit the risk of obtaining physically unrealizable solutions, a DC approximation of AC power equations is commonly used to reduce the computational complexity \cite{fernandez2016optimal,7317588}. AC OPF models are solved with simulations \cite{6736137}, local solvers, and heuristic methods \cite{crossland2014planning,7052416}, which cannot guarantee the global optimality of the proposed solutions. Recently, convex relaxations of the OPF problem have drawn research interest since they can produce globally optimal solutions. Mainly semidefinite programming and second-order cone programming (SOCP) have been widely studied \cite{zohrizadeh2020survey}. In the context of renewable and storage integration, convex relaxations become harder to solve since the decisions require solving a mixed-integer multi-period OPF (MOPF). The studies by \cite{7389434} and \cite{grover2018optimal} are two of the few papers that propose an exact SOCP relaxation for energy storage optimization. However, the exactness of SOCP is conditioned on certain settings and valid for only radial distribution systems. Effective convex programming approaches for meshed transmission networks need further development.

Hydrogen storage, as opposed to other alternative storage systems, interacts with the external hydrogen market and provides opportunities for selling hydrogen. For example, \cite{ghofrani2013framework} shows that arbitrage revenues alone cannot justify the investment cost of storage under some settings. However, the cost may be justified by considering other storage-related benefits, such as profit from selling hydrogen. Therefore, consideration of the hydrogen market has the potential to change decision dynamics in the context of storage location and sizing, which is yet unaddressed in the literature. 

In this paper, we propose a stochastic optimization model for jointly deciding on the location and sizing of renewables and hydrogen storage based on multi-period AC OPF problems. We propose a solution approach based on SOCP to provide globally optimal solutions for the resulting model. We create a representative test case that involves scenarios based on real data sets to represent the stochasticity of the network load, renewable energy output, and electricity generation price. To the best of our knowledge, this is the first study to provide a systematic optimization method to decide on optimal locations and sizings of renewables and hydrogen storage while considering the various dynamics of the underlying operational level problem including AC power flow equations, the stochasticity of operational parameters, and integration of the hydrogen market. Specifically, the following strategic level questions can be answered: (1) Can operational cost savings compensate for the high capital costs of renewables and hydrogen storage? (2) How should investment budgets be allocated between renewables and hydrogen storage? (3) Which locations and sizing are preferable for renewables and hydrogen storage? The main contributions can be summarised as follows:
\begin{itemize}
    \item We propose a new stochastic optimization model for joint renewable and hydrogen storage location and sizing into power networks based on multi-period AC OPF problems. In addition, our model captures interaction with the hydrogen market.
    \item We develop a systematic solution approach based on SOCP within a Benders decomposition framework to provide globally optimal solutions. Our approach offers global optimality guarantees with very small optimality gaps.
    \item  On a representative test case, we show it is crucial to consider the joint optimization of renewables and hydrogen storage as it results in significant operational cost savings compared to the case where we only include renewables. Moreover, by comparing against DC approximations, we show the importance of including AC power equations as it changes location and sizing decisions and thereby reduces operational costs.
    \item Our optimization framework allows us to answer relevant strategic-level questions. Namely, we show that a functioning hydrogen market can change decision dynamics. In addition, we investigate the effects of renewable curtailment cost, emission tax price, and conversion efficiencies by means of a sensitivity analysis. Our findings are useful for decision-makers in integrating renewables and hydrogen storage in power networks.
\end{itemize}

The remainder of this paper is organized as follows. Section~\ref{sec:model} introduces the optimization model with its SOCP relaxation. Section~\ref{sec:solution method} introduces our solution approach for this optimization problem. Section~\ref{sec:model input} describes the model input and introduces the input data used in the model formulation. Section~\ref{sec:computational experiments} presents the computational results. Section~\ref{sec:conclusions} presents the concluding remarks.

\section{Model}\label{sec:model}

This section presents our mathematical programming formulation and its mixed-integer SOCP (MISOCP) relaxation. The system consists of three main components: a power network, renewable energy sources, and hydrogen storage. The interaction of these components is coordinated by a central network operator responsible for investing in renewable energy sources, investing in hydrogen storage systems, and planning the daily network operations. We consider that the investment decisions are made once to be operational during its lifetime. To simulate the operation of the resulting power system after investment decisions are made, representative days are used to characterize the daily network planning. The goal is to minimize the expected daily operational cost for a given investment budget. We model the joint optimization problem of the network operator as a two-stage stochastic mixed-integer non-linear programming (MINLP) model.

The power network is denoted by $\mathcal{N} = (\cB, \cL)$, where $\cB$ denotes the set of buses and $\cL$ denotes the set of transmission lines. Let $\delta(i)$ denote the set of neighbors for bus $i \in \mathcal{B}$ and let $\cG \subseteq \cB$ denote the set of conventional (i.e., non-renewable) generators. 

The first-stage decision comprises the location and sizing of renewable energy sources (e.g., wind turbines) and hydrogen storage subject to a given investment budget. Hydrogen storage consists of an electrolyzer to convert renewable power into hydrogen, a storage unit to store hydrogen, and a fuel cell to convert the hydrogen back into power. We assume that these components are installed together, and the storage unit and fuel cell capacity are in line with the size of the electrolyzer. The second-stage decisions take place after the uncertainty of renewable output, electricity demand, and electricity generation prices are revealed. It entails planning daily network operations by solving the AC MOPF over a finite time horizon $\cT = \{1,...,T\}$ subject to a given set of scenarios (i.e., representative days) $\omega \in \Omega$. The network operator can decrease the cost of daily network operations by exploiting the energy arbitrage by storing electricity when prices are low and feeding back electricity to the network when prices are high. Moreover, the profits can be boosted by selling hydrogen to the external market.

In what follows, we first detail the two-stage stochastic MINLP model that has a non-convex feasible region due to the AC power equations. Afterward, we present its convex relaxation based on SOCP.

\subsection{MINLP Formulation}

The objective function~\eqref{eq:objective} of the MINLP formulation minimizes the expected operational cost, which consists of five parts: the cost function of production from conventional generators ($h(\cdot)$), a penalty term representing the emission cost associated with conventional generators ($C^E$), the cost of curtailing excess production from renewables ($C^C$), the cost of unsupplied loads ($C^U$), and profit obtained from selling hydrogen to external market ($C^S$).

\begin{equation}\begin{split}\label{eq:objective}
\allowdisplaybreaks
\min  \hspace{0.5em}  & \sum_{\omega\in\Omega}  \rho_{\omega} \biggl[\sum_{t\in\cT} \left(  \sum_{i\in\cG} \left( h(p_{it}^g(\omega)) + C^E p_{it}^g(\omega) \right) \right)    \\ &   + \left(  \sum_{i\in\cB} \left(  C^C l_{it}(\omega) + C^U u_{it}(\omega)  - C^S h_{it}(\omega) \right)\right) \biggr]. 
  \end{split}\end{equation}
The system is subject to the following constraints: 
\subsubsection{Investment Constraints}
 We denote the location decisions for renewable and hydrogen storage with $r_i^{\textsc{b}}$ and $h_i^{\textsc{b}}$, respectively, equaling $1$ if a new source is located to bus $i \in \mathcal{B}$, and 0 otherwise. We determine the corresponding power ratings with continuous decision variables $r_i^{\textsc{r}}$ and $h_i^{\textsc{r}}$.  
\begin{subequations}\label{Cset:investments}
\allowdisplaybreaks\begin{align}
&\sum_{i \in \cB} (C^R r_{i}^{\textsc{r}} + C^H h_{i}^{\textsc{r}}) \leq \overline{B} \ &&  \label{eq:budget constraint} \\
&  h_{i}^{\textsc{b}} \le r_{i}^{\textsc{b}} &i&\in \cB \label{eq:green hydrogen constraint} \\
&\underline{r_{i}^{\textsc{r}}} r_{i}^{\textsc{b}} \le r_{i}^{\textsc{r}} \le \overline{r_{i}^{\textsc{r}}} r_{i}^{\textsc{b}}  &i&\in \cB \label{eq:renewable capacity limits} \\
&\underline{h_{i}^{\textsc{r}}} h_{i}^{\textsc{b}} \le h_{i}^{\textsc{r}} \le  \overline{h_{i}^{\textsc{r}}} h_{i}^{\textsc{b}}  &i&\in \cB. \label{eq:hydrogen storage capacity limits} 
\end{align}
\end{subequations}
Constraint~\eqref{eq:budget constraint} limits the total investments in renewables and hydrogen storage by a certain budget ($\overline B$). Constraint~\eqref{eq:green hydrogen constraint} limits placing hydrogen storage to a node with a renewable energy source. Constraints~\eqref{eq:renewable capacity limits} and \eqref{eq:hydrogen storage capacity limits} ensure that the power ratings of renewables and hydrogen storage are within the prespecified ranges, respectively.
 
\subsubsection{Operational Storage-related Constraints}
For each bus $i \in \cB$, time $t \in \cT$, and scenario $\omega \in \Omega$:
\begin{subequations}\label{Cset:hydrogen}
\allowdisplaybreaks\begin{align}
    & s_{it}(\omega)+\eta_{g} (p2g)_{it}(\omega)-(g2p)_{it}(\omega)-h_{it}(\omega)=s_{i(t+1)}(\omega)  \label{eq:stock constraint}\\
    &  s_{i0}(\omega) = I_{i}(\omega)  \label{eq:initial stock}\\
    &  s_{it}(\omega) \leq s^{\text{max}} h_{i}^{\textsc{r}}   \label{eq:stock limit}\\
    &  (p2g)_{it}(\omega) \leq h_{i}^{\textsc{r}} \label{eq:p2g limit}\\
    &  (p2g)_{it}(\omega) \leq r_{it}(\omega) r_{i}^{\textsc{r}}  \label{eq:p2g limit2}\\
    &  (g2p)_{itk} \le f^{\text{rate}} h_{i}^{\textsc{r}}. \label{eq:g2p limit}
    \end{align}
\end{subequations}
Constraint~\eqref{eq:stock constraint} controls the hydrogen level between consecutive periods by considering the amount of power-to-gas, gas-to-power conversions, and the selling of hydrogen. Constraint~\eqref{eq:initial stock} sets the hydrogen storage's initial state of charge. Constraint \eqref{eq:stock limit} ensures that the storage capacity is not exceeded. Constraints~\eqref{eq:p2g limit} and \eqref{eq:g2p limit} limit the power-to-gas conversion and gas-to-power conversions, respectively. Constraint~\eqref{eq:p2g limit2} allows only renewable power to be converted into green hydrogen.

\subsubsection{Node Balance Constraints}  
For each bus $i \in \cB$, time $t \in \cT$, and scenario $\omega \in \Omega$:
    \begin{subequations}\label{Cset:node-balance}
\allowdisplaybreaks\begin{align}
    &\hspace{0.5em} p_{it}^g(\omega) - p_{it}^d(\omega) +r_{it}(\omega)r_i^{\textsc{r}} -l_{it}(\omega)-(p2g)_{it}(\omega) \label{eq:active balance}  \\
    &\hspace{0.4cm}+(g2p)_{it}(\omega) \eta_{p} + u_{it}(\omega)=  \ g_{ii} |V_{it}(\omega)|^2 +  \sum_{j\in\delta(i)} p_{ijt}(\omega) \   \nonumber    \\
    &\hspace{0.5em} q_{it}^g(\omega)  - q_{it}^d(\omega)  = -b_{ii} |V_{it}(\omega) |^2 + \sum_{j\in\delta(i)} q_{ijt}(\omega).    \label{eq:reactive balance}
\end{align}
\end{subequations}
Constraint~\eqref{eq:active balance} ensures active power flow balance at bus $i$ while considering uncertain network load, uncertain renewable output and curtailments, power-to-gas and gas-to-power conversions, and unsupplied load. Constraint~\eqref{eq:reactive balance} ensures reactive power flow balance at bus $i$. 
\subsubsection{Flow Constraints} For each line $ (i,j) \in \cL$, time $t \in \cT$, and scenario $\omega \in \Omega$:
    \begin{subequations}\label{Cset:Flow}
\allowdisplaybreaks\begin{align}
    &  p_{ijt}(\omega) = G_{ij}|V_{it}(\omega)|^2 + |V_{it}(\omega)| |V_{jt}(\omega)|  \label{eq:active flow}   \\
    & \hspace{0.5cm} \times [G_{ij} \cos(\theta_{it}(\omega) - \theta_{jt}(\omega))-B_{ij} \sin(\theta_{it}(\omega) - \theta_{jt}(\omega)) ] \nonumber \\
    &  q_{ijt}(\omega) = - B_{ij}|V_{it}(\omega)|^2  - |V_{it}(\omega)| |V_{jt(\omega)}| \label{eq:reactive flow}\\ 
    & \hspace{0.5cm} \times [B_{ij}\cos(\theta_{it}(\omega) - \theta_{jt}(\omega)) +G_{ij} \sin(\theta_{it}(\omega) - \theta_{jt}(\omega)) ]. \nonumber
\end{align}
\end{subequations}
Constraints~\eqref{eq:active flow} and \eqref{eq:reactive flow} represent the active and reactive power flow, respectively.
\subsubsection{Network Operational Limits} For each time $t \in \cT$, and scenario $\omega \in \Omega$:
     \begin{subequations} \label{Cset:operational}
\allowdisplaybreaks\begin{align}  
    & \underline V_i^2 \le |V_{it}(\omega)| \le \overline  V_i^2  \ & i& \in \cB \label{eq:voltage limit}\\
    & \underline{p}_i  \le p_{it}^g(\omega) \le \overline{p}_i &g& \in \cG \label{eq:active gen}\\
    &  \underline{q}_i  \le q_{it}^g(\omega) \le \overline{q}_i &g& \in \cG  \label{eq:reactive gen}\\
    &  -R_i^{\text{down}} \leq p_{it+1}^g(\omega) - p_{it}^g(\omega)    \leq R_i^{\text{up}}&i& \in \cB \label{eq:ramp limits} \\
    &  p_{ijt}(\omega)^2+q_{ijt}(\omega)^2  \le \overline S_{ij}^2 &(&i,j) \in \cL   \label{eq:flow limit}\\
    &  |\theta_{it}(\omega) - \theta_{jt}(\omega)| \le \overline\theta_{ij} &(&i,j) \in \cL. \label{eq:angle limit}
\end{align}
\end{subequations}
 Constraint~\eqref{eq:voltage limit} enforce bus voltage magnitude to stay within acceptable limits of lower and upper bounds. Constraints~\eqref{eq:active gen} and \eqref{eq:reactive gen} limit the active and reactive power outputs of generator $i$. We set $\underline{p}_i = \overline{p}_i =\underline{q}_i = \overline{q}_i = 0$ for $ i \in \mathcal{B} \setminus  \mathcal{G}$.  Constraint~\eqref{eq:angle limit} sets the ramp down and ramp up limits of generator $i$. Constraints~\eqref{eq:flow limit} and \eqref{eq:angle limit} limit the transmission capacity and the phase angle of line $(i,j)$ as a function of the maximum allowable flow and phase angle bound, respectively.

The MINLP formulation is obtained as $M_O$:\{\eqref{eq:objective} : \eqref{Cset:investments}--\eqref{Cset:operational}\}

\subsection{An Alternative Formulation}
In order to obtain an SOCP-based relaxation for the MINLP model problem, we first provide an alternative formulation motivated by \cite{kocuk2016strong,kayacik2021promise}. We define the following decision variables:
\begin{itemize}
\item 
For each bus $i \in \cB$, time $t \in \cT$, and scenario $\omega \in \Omega$,
\begin{itemize}
\item $c_{iit}(\omega) := |V_{it}(\omega)|^2$.
\end{itemize}
\item
For each line $(i,j) \in \cL$, time $t \in \cT$, scenario $\omega \in \Omega$ ,
\begin{itemize}
\item
$c_{ijt}(\omega) := \ |V_{it}(\omega)||V_{jt}(\omega)| \cos(\theta_{it(\omega)} - \theta_{jt}(\omega))$
\item
$s_{ijt}(\omega) := -|V_{it}(\omega)||V_{jt}(\omega)| \sin(\theta_{it(\omega)} - \theta_{jt(\omega)})$.
\end{itemize}
\end{itemize}

We substitute the new variables in Constraints~\eqref{Cset:node-balance}, \eqref{Cset:Flow}, and \eqref{eq:voltage minlp} and linearize them as follows. For each time $t \in \cT$ and scenario $\omega \in \Omega $:
\begin{subequations}\label{Cset:updatedSOCP}
\allowdisplaybreaks\begin{align}
 &p_{it}^g(\omega) - p_{it}^d(\omega) +r_{it}(\omega)r_i^{\textsc{r}} -l_{it}(\omega)-(p2g)_{it}(\omega) \label{eq:active balance minlp}  \\
 &+(g2p)_{it}(\omega)\eta_p + u_{it}(\omega) = g_{ii} c_{iit}(\omega) +  \sum_{j\in\delta(i)} p_{ijt}(\omega)  \hspace{0.1cm}  i \in \cB \nonumber \\
    &\ q_{it}^g(\omega) - q_{it}^d(\omega)= - {b_{ii}} {c_{iit}(\omega)}  + \sum_{j\in\delta(i)} q_{ijt}(\omega)  \hspace{0.4cm} i \in \cB \label{eq:reactive balance minlp_}  \\
    & p_{ijt}(\omega) = G_{ij} c_{iit}(\omega)  + G_{ij}c_{ijt}(\omega)    - B_{ij}s_{ijt}(\omega)  \hspace{0.2cm}   (i,j) \in \cL \label{eq: active flow minlp_}\\
    & q_{ijt}(\omega) = -B_{ij} c_{iit}(\omega)  -  B_{ij}c_{ijt}(\omega)  - G_{ij} s_{ijt}(\omega) \hspace{0.1cm}   (i,j) \in \cL \label{eq: reactive flow minlp_} \\
     & \underline V_i^2 \le c_{iit}(\omega) \le \overline  V_i^2 \hspace{4.2cm}  i \in \cB. \label{eq:voltage minlp} 
\end{align}
\end{subequations}

To preserve the trigonometric relation between the new variables $c_{iit}(\omega), c_{ijt}(\omega), s_{ijt}(\omega)$, we need additional non-convex constraints. These, so-called consistency constraints are defined for each line $(i,j) \in \cL$, time $t \in \cT$, and scenario $\omega \in \Omega$ as follows:
\begin{subequations}\label{Cset:consistency}
\allowdisplaybreaks\begin{align}
 & c_{ijt}(\omega)^2+s_{ijt}(\omega)^2  = c_{iit}(\omega) c_{jjt}(\omega)  \label{eq:consistency1}    \\
&\theta_{jt}(\omega)- \theta_{it}(\omega) = \atantwo(s_{ijt}(\omega),c_{ijt}(\omega)). \label{eq:consistency2}
\end{align}
\end{subequations}

Subsequently, an alternative exact formulation to $M_O$ is obtained as: \{\eqref{eq:objective}: \eqref{Cset:investments}, \eqref{Cset:hydrogen}, \eqref{eq:active gen}--\eqref{eq:angle limit}, \eqref{Cset:updatedSOCP}, \eqref{Cset:consistency}\} 

\subsection{MISOCP Relaxation}
We convexify the consistency constraints by eliminating Constraint~\eqref{eq:consistency2} and relaxing Constraint~\eqref{eq:consistency1} as follows:
\begin{equation}\label{eq:consistency relaxed}
    c_{ijt}(\omega)^2+s_{ijt}(\omega)^2  \le c_{iit}(\omega) c_{jjt}(\omega).   
\end{equation}
The MISOCP relaxation of the proposed formulation is obtained as $M_{R}$ : \{\eqref{eq:objective}: \eqref{Cset:investments}, \eqref{Cset:hydrogen}, \eqref{eq:active gen}--\eqref{eq:angle limit}, \eqref{Cset:updatedSOCP}, \eqref{eq:consistency relaxed}\}. \\

\section{Solution Method}\label{sec:solution method}
We propose a systematic solution method based on MISOCP. The original problem $M_O$ is challenging to solve with standard local solvers; even if solved, a locally optimal solution can be obtained. Therefore, we use the MISOCP relaxation $M_{R}$ to aim for globally optimal solutions to $M_O$. If the convex relaxation is exact, it guarantees global optimality to the original problem. Although the SOCP relaxation of OPF is rarely exact in practice, we can still exploit it in two aspects: First, it provides a lower bound (LB) for the optimal value of the original problem $M_O$. Second, we utilize the optimal solution of the relaxation to guide a local solver to obtain a feasible solution, hence, an upper bound (UB), for the original problem $M_O$. In this way, we obtain lower and upper bounds to $M_O$, from which we can compute a quality measure for global optimality.

The MISOCP relaxation $M_{R}$ is a stochastic multi-period mixed-integer programming model, and it is hard to solve using standard solvers (e.g., Gurobi) for increasing instance size. Therefore, we propose Benders decomposition to solve $M_{R}$. We first separate the model into a master problem ($MP$) and $|\Omega|$ subproblems ($SP_{\omega}$). In the master problem, we make the location and sizing decisions subject to the investment constraints (Constraints~\eqref{Cset:investments}). In each subproblem, we solve the SOCP relaxation of the multi-period OPF for a fixed scenario $\omega \in \Omega$ ($SP_{\omega}$ : \{\eqref{eq:objective} : \eqref{Cset:hydrogen}, \eqref{eq:active gen}--\eqref{eq:angle limit}, \eqref{Cset:updatedSOCP}, \eqref{eq:consistency relaxed}\}). Note that Benders decomposition converges to an optimal solution if the subproblems are convex \cite{geoffrion1972generalized}. 

Algorithm~\ref{SolutionAlgorithm} details our solution method. In the first step, we solve master problem $MP$ to obtain an initial solution set of investment decisions $P^* = \{r_i^{\textsc{b}*},r_i^{\textsc{r}*},h_i^{\textsc{b}*},h_i^{\textsc{r}*}\}$ and a Benders lower bound (BLB). Then, we solve each subproblem $SP_{\omega}$ to obtain Benders upper bound (BUB) and generate optimality cuts ($\Phi$). Then, we include the optimality cuts in $MP$ and solve the resulting problem to update BLB. We repeatedly solve the $SP_{\omega}$'s and $MP$ until the Benders optimality gap is smaller than $\epsilon$. In the second step, we set our global lower bound LB equal to BLB. We then fix the investment decisions $P$ in the $M_O$ to obtain UB from the remaining non-linear program (NLP). Since the investment decisions are fixed in $M_O$, the remaining problem becomes an NLP that can be decomposed into $|\Omega|$ subproblems as $M_O{\omega}$ for each scenario $\omega \in \Omega$. We solve these subproblems $M_O{\omega}$ to obtain UB. Lastly, we calculate the global optimality gap. 

\begin{algorithm}
\caption{Solution Approach}
\label{SolutionAlgorithm}
\begin{algorithmic}[1]
\STATE Set: LB = BLB = -$\infty$, UB = BUB = $\infty$
\STATE \underline{\textbf{Step 1:}}
\STATE Solve $MP$ to obtain  $P^* = \{r_i^{\textsc{b}*},r_i^{\textsc{r}*},h_i^{\textsc{b}*},h_i^{\textsc{r}*}\}$
\STATE  Set BLB $\leftarrow z(MP)$
\WHILE {$(1-\text{BLB}/\text{BUB})<\epsilon $}
 \FORALL{$\omega \in \Omega\ $}
 \STATE  $z(SP_{\omega}) \leftarrow$ Solve $SP_{\omega}$ subject to $P^*$.
 \ENDFOR
 \STATE BUB $\leftarrow \min\left({\sum_{\omega \in \Omega}z(SP_{\omega}),\text{BUB}}\right)$.
 \IF{BUB $\geq \sum_{\omega \in \Omega}z(SP_{\omega})$}
 \STATE $P \leftarrow P^*$
 \ENDIF
 \STATE Generate optimality cuts $\Phi$
 \STATE Solve MP with $\Phi$ to obtain $P^*$
\STATE  Set BLB $\leftarrow z(MP)$
 \ENDWHILE 
 \STATE \underline{\textbf{Step 2:}} 
 \STATE Set LB to BLB
 \STATE Fix $P$ and decompose $M_{O}$ into $|\Omega|$ subproblems as $M_{O_{\omega}}$
\FORALL{$\omega \in \Omega\ $}
\STATE $z(M_{O_{\omega}}) \leftarrow $ Solve $M_{O_{\omega}}$
\ENDFOR
\STATE UB $\leftarrow \sum_{\omega \in \Omega} z(M_{O_{\omega}})$
\STATE Compute global optimality gap as $100\times(1-\text{LB}/\text{UB})$ 
    \end{algorithmic}
\end{algorithm}

\section{Model Input}\label{sec:model input}
Our model input is based on discussions with the stakeholders in the energy sector within the HEAVENN Program in the Northern Netherlands, where Europe's first hydrogen valley is being built \cite{Heavenn}. Our model draws on data on the power network dynamics, renewable energy production, and hydrogen demand. We detail the related data in this section, while parameters for sensitivity analysis are in Section~\ref{sec:computational experiments}. 

For the MOPF dynamics, we consider a 24-hour time horizon of 1-hour periods, i.e., $|\cT| = 24$, from 00:00 to 00:00 the following day. We create daily scenarios to specify realized electricity demand, electricity price, and renewable energy supply for each of the 24-hour periods. Each scenario represents a typical day of each season of the year 2021, resulting in four representative days throughout the year. 
 
\subsection{OPF instance}
As the actual grid data of the Netherlands is confidential, we test our algorithm on the well-established OPF instance IEEE30 from the Power Grid Library (PGLIB-OPF) \cite{birchfield2016grid}, which includes the network structure and parameters for a single period. To make it compatible with the multi-period formulation, we adjust the relevant parameters that vary on an hourly basis (e.g., electricity demand and price) and keep other parameters fixed.

We obtain electricity demand and day-ahead electricity prices data from the European Network of Transmission System Operators for Electricity (ENTSOE) \cite{ENTSOE2022}. We calculate the hourly averages of electricity demand for each season associated with the representative day. After normalizing the hourly averages by their maximum, we multiply the network's active and reactive power load with the corresponding normalized values. We ended up with an average of $195$ kWh hourly active power load ranging from $150$ to $283$kWh. See Algorithm~\ref{Algo:ElectricityDemand} for the details.

\begin{algorithm}
\caption{Hourly electricity demand }
\label{Algo:ElectricityDemand}
\textbf{Input:} From ENTSOE : Set of hourly average electricity demand for each scenario 
$  \mathcal{D} =  \{d_{t}(\omega) : t \in \mathcal{T}, \omega \in \Omega \}$ \ From OPF Data : $p_i^d, q_i^d$  \\
\textbf{Output:} Hourly power load values $p_{it}^d(\omega),  q_{it}^d(\omega) $
\begin{algorithmic}[1]
\STATE $\texttt{max\_demand}= \max( \mathcal{D})$
\FORALL{$t \in \mathcal{T}, \omega \in \Omega$}
\STATE $\overline{d}_{t}(\omega)= \frac{d_{t}(\omega)}{\texttt{max\_demand}} $
\FORALL{$i \in \mathcal{B}$}
\STATE $p_{it}^d(\omega) = p_i^d \times  \overline{d}_{t}(\omega)$ 
\STATE $q_{it}^d(\omega) = q_i^d \times  \overline{d}_{t}(\omega)$
\ENDFOR
\ENDFOR
\end{algorithmic}
\end{algorithm}

For hourly electricity generation prices, we multiply the hourly average day-ahead prices for each season by the normalized costs of generators over the whole network. Details are outlined in Algorithm~\ref{Algo:ElectricityPrice}. We attain an average of $80$ €/MWh hourly electricity generation price ranging from $20$ to $220$ €/MWh.

\begin{algorithm}
\caption{Hourly electricity generation cost}
\label{Algo:ElectricityPrice}
\textbf{Input:} From ENTSOE: Daily day-ahead prices for each scenario  
$  \mathcal{P} =  \{p_{t}(\omega) : t \in \mathcal{T}, \omega \in \Omega \}$, From OPF Data : Generation cost for each generator $  \mathcal{C} =  \{c_{g} : g \in \mathcal{G} \}$   \\
\textbf{Output:} Generation cost for each generator, for each time, and for each scenario $\mathcal{H}=\{h_{gt}(\omega): g\in\mathcal{G}, t \in \mathcal{T}, \omega \in \Omega \}$
\begin{algorithmic}[1]
\STATE $\texttt{average\_cost}= \text{average}( \mathcal{C})$
\FORALL{$t \in \mathcal{T}, k \in \mathcal{K}$}
\STATE $h_{gt}(\omega)= \frac{c_g}{\texttt{average\_cost}} \times p_{t}(\omega) $
\ENDFOR
\end{algorithmic}
\end{algorithm}

We set the cost of unsupplied demand ($C^U$) to $3000$ €/MWh based on \cite{7509683}. 
\subsection{Renewable Data}
We obtain hourly wind speed data from  the Koninlijk Nederlands Meteorologisch Instituut (KNMI) \cite{KNMI2022}. We consider a wind turbine with the specifications of a Vestas V20 with a $4.5$ m/s cut-in wind speed ($v_{\mathrm{ci}}$), a $13.0$ m/s rated wind speed($v_{\mathrm{r}}$), and $25.0$ m/s as the cut-out wind speed ($v_{\mathrm{co}}$) \cite{vestasV20}. Given the hourly wind speed ($v_{it}(\omega)$) from KNMI data, we calculate the hourly wind power factor ($r_{it}(\omega)$) of a wind turbine as in \cite{atwa2011probabilistic}. 

$ r_{it}(\omega)=\left\{\begin{array}{ll}0, & 0 \leq v_{it}(\omega)  \leq v_{\mathrm{ci}} \\  \frac{\left(v_{it}(\omega) -v_{\mathrm{c}}\right)}{\left(v_{\mathrm{r}}-v_{\mathrm{ci}}\right)}, & v_{\mathrm{ci}} \leq v_{it}(\omega)  \leq v_{\mathrm{r}} \\ 1, & v_{\mathrm{r}} \leq v_{it}(\omega)  \leq v_{\mathrm{co}} \\ 0, & v_{\mathrm{co}} \leq v_{it}(\omega) \end{array}\right.$

We assume that a wind turbine with a minimum size of $100$ kW can be installed and accordingly set minimum power rating $\underline{r_i^{\textsc{r}}}$ to $100$ kW. We consider the capital cost of a wind turbine with a lifetime of $20$-$30$ years as $1.2$ M€/MW and set $C^R$ to $1.2$ M€/MW \cite{IRENA2021}. 

\subsection{Hydrogen Data}
The capital cost and lifetime of the electrolyzer, fuel cell, and hydrogen storage tank are provided in Table~\ref{tab:Hydrogen Data}. We consider average capital costs values corresponding to $1.05, 18.75\times 10^{-6}$, and $1.05$ M€/MW for electrolyzer, storage, and fuel cell, respectively. We adjust the remaining hydrogen-related parameters based on \cite{li2021value}. We assume that the power rating of the electrolyzer is at least $30 \% $ of the minimum power rating of wind turbine, and set $\underline{h_i^{\textsc{r}}}$ to $30$ kW, which corresponds to $30 \% $ of the minimum power rating of renewables $\underline{r_i^{\textsc{r}}}$. Based on a setting in which the fuel cell capacity is half of the installed electrolyzer capacity, and the storage tank can store $20$ hours of full electrolyzer output, we set $f^{\text{rate}}$ to $0.5$ and $s^{\text{max}}$ to $20$. Accordingly, a storage unit with 1 MW of electrolyzer costs approximately $1.6$ €/MW, and  $C^H$ corresponds to 1.6 €/MW. We set the conversion efficiencies of electrolyzer $\eta_g$ and fuel cell  $\eta_p$ to $0.7$ and $0.5$, respectively.
\begin{table}[h]
\centering
    \caption{ Hydrogen Data }
    \label{tab:Hydrogen Data}
\begin{tabular}{|l|l|l|}
\hline
Technology & Capital cost & Lifetime \\
\hline
Electrolyzer \cite{le2017relevance} 
& 0.7 - 1.4 (M€/MW) & 20-30y   \\
Storage  \cite{le2017relevance} &  (6.75- 30.75) $\times 10^{-6}$ (M€ /MW)  & 50y  \\
Fuel cell \cite{schmidt2017future} & 0.7 - 1.4 (M€ /MW) & 20-30y  \\

\hline
\end{tabular}
\end{table}

\section{Computational Experiments}\label{sec:computational experiments}
We present the results in two parts. First, we show the trade-off between investment and operational costs. Second, we present the corresponding optimal location and power ratings. All computational experiments are carried out on an Intel Xeon E5 2680v3 CPU with a 2.5 GHz processor and 32 GB RAM. Implementation is coded in Python with Gurobi 9.1.0 and IPOPT for solving the MISOCP relaxation and the NLP models, respectively. We note an average of 3\% global optimality gap for all the settings.

\subsection{Trade-off Between Investment and Operational Costs}
In this section, we vary the several parameters that can be influenced by economic, technical, and regulatory policies to analyze their effect on investment decisions. Figures~\ref{fig:HydrogenPriceFigures}--\ref{fig:EfficiencyFigures2} display results for varying hydrogen selling prices, curtailment costs, emission tax prices, and conversion efficiencies. In each figure, the graph on the left shows the trade-off between investment budget and operational cost, while the graph on the right shows the percentage of the investment budget allocated to hydrogen storage. Unless otherwise stated, the hydrogen selling price is set to $2$ €/kg, curtailment cost to $40$ €/MWh, and emission tax price to $30$ €/ton \ce{\ce{CO2} }  (which corresponds to the setting with red dashed lines in each figure). In each figure, we vary one parameter between the ranges specified in Table~\ref{tab:SensitivityParameters}, which are deemed relevant values based on discussions within the HEAVENN Program. The emission tax price is derived from the proposed emission tax price by the Dutch government, which is €$30$ per ton of \ce{\ce{CO2} }  for 2021 and €$125$ for 2030 \cite{Klimaatakkoord2019}.
\begin{table}[h]
\centering
    \caption{Sensitivity analysis parameters}
    \label{tab:SensitivityParameters}
\begin{tabular}{|p{0.65\linewidth}|l|}
\hline
Parameter & Range  \\
\hline
Green hydrogen selling price ($C^S)$ &  0-6 €/kg \\
Curtailment cost ($C^C$) & 0-120 €/MWh \\
Emission tax price ($C^E$) &  0-125€/ton \ce{\ce{CO2} }  \\
Power-to-gas, gas-to-power efficiencies ($\eta_g,\eta_p$): &  0.7-1, 0.5-1 \\
\hline
\end{tabular}
\end{table}
We obtain the graphs by solving the model under different investment budgets ranging from 0 to 1M€. In the graphs, we also provide the investment budgets scaled to a daily basis to better reflect the overall daily expense of the power system throughout the lifespan of green technologies. Accordingly,  we use the following equation: 
$$
D_c= C \frac{\delta \cdot(1+\delta)^{\gamma}}{(1+\delta)^{\gamma}-1} \cdot \frac{1}{N_{\text {year }}}
$$
where $D_c$ is the daily capital cost, $C$ is the capital cost, $\delta$ is the annual discount rate, $\gamma$ is the lifetime, and  ${N_{\text {year }}}$ is the number of days in a year. Based on Table~\ref{tab:Hydrogen Data} and \cite{IRENA2021}, we assume that a wind turbine and a hydrogen storage unit have a lifetime of $\gamma=25$ years and an annual discount rate of $\delta= 5\%$. While the actual investment budget ranges from zero to 1M€, corresponding daily scaled values range from 0 to 200€. 

\subsubsection{Effect of Hydrogen Market}
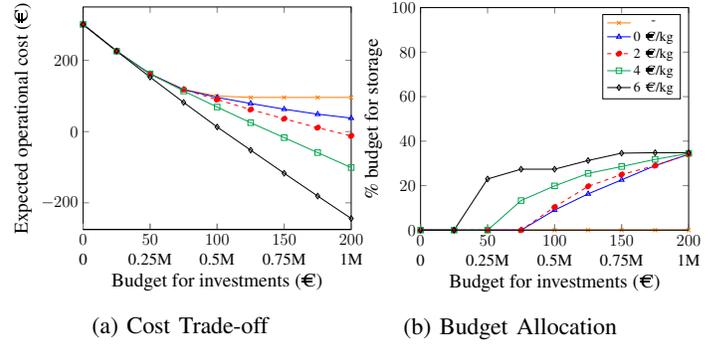
\begin{figure}[h]
\begin{subfigure}{0.25\textwidth}
  \centering
  \begin{tikzpicture}[scale=0.52]
\begin{axis}[
    ylabel={\Large Expected operational cost (€)},
    xlabel={\Large Budget for investments (€)},
    y tick label style={font=\large},
    x tick label style={font=\large},
    xmin=0, xmax=200,
    xlabel shift={14pt},
    ymin=-275, ymax=350,
    xtick={0,50,100,150,200},
    legend pos=south west,
    grid style=dashed,
]

\addplot[
    color=orange,
    mark=x,
    ]
    coordinates {
  (0,301)
  (25,226)
  (50,162)
  (75,118)
  (100,100)
  (125,96)
  (150,96)
  (175,96)
  (200,96)
    };
    
\addplot[
    color=blue,
    mark=triangle,
    ]
    coordinates {
  (0,301)
  (25,226)
  (50,162)
  (75,118)
  (100,96)
  (125,79)
  (150,63)
  (175,49)
  (200,38)
    };
    
\addplot[
    color=red,
    mark=*,
    style=dashed,
    ]
    coordinates {
  (0,301)
  (25,226)
  (50,162)
  (75,118)
  (100,90)
  (125,62)
  (150,36)
  (175,11)
  (200,-12)
    };
    
\addplot[
    color=Green,
    mark=square,
    ]
    coordinates {
  (0,301)
  (25,226)
  (50,162)
  (75,114)
  (100,69)
  (125,25)
  (150,-17)
  (175,-59)
  (200,-101)
    };

\addplot[
    color=black,
    mark=diamond,
    ]
    coordinates {
  (0,301)
  (25,226)
  (50,153)
  (75,82)
  (100,13)
  (125,-52)
  (150,-117)
  (175,-181)
  (200,-244)
    };

    
    
    
    

    \end{axis}
    \begin{axis}[
        xmin=0,xmax=1,
        xtick={0,0.25,0.5,0.75,1},
        xticklabels={0,0.25M,0.5M,0.75M,1M},
        hide y axis,
        axis x line*=none,
        x tick label style={font=\large},
        ymin=0, ymax=32000,
        x label style={yshift=-0.5cm},
        x tick label style={yshift=-0.5cm}
        ]
    \end{axis} 
\end{tikzpicture}
\caption{Cost Trade-off} \label{fig:HydrogenPriceTrade-off} 
\end{subfigure}
\begin{subfigure}{.22\textwidth}
  \centering
  \begin{tikzpicture}[scale=0.52, transform shape]
\begin{axis}[
    ylabel={\Large \% budget for storage},
    xlabel={\Large Budget for investments (€)},
    y tick label style={font=\large},
    x tick label style={font=\large},
    ymin=0, ymax=100,
    xmin=0, xmax=200,
    xlabel shift={14pt},
    xtick={0,50,100,150,200},
    legend pos=north east,
    grid style=dashed,
]

\addplot[
    color=orange,
    mark=x,
    ]
    coordinates {
  (0,0)
  (25,0)
  (50,0)
  (75,0)
  (100,0)
  (125,0)
  (150,0)
  (175,0)
  (200,0)
    };
    
\addplot[
    color=blue,
    mark=triangle,
    ]
    coordinates {
   (0,0)
  (25,0)
  (50,0)
  (75,0)
  (100,9)
  (125,16.3)
  (150,22.6)
  (175,29)
  (200,34.1)
    };
    
\addplot[
    color=red,
    mark=*,
    style=dashed,
    ]
    coordinates {
  (0,0)
  (25,0)
  (50,0)
  (75,0)
  (100,10.4)
  (125,19.7)
  (150,25)
  (175,29)
  (200,34.1)
    };
    
\addplot[
    color=Green,
    mark=square,
    ]
    coordinates {
  (0,0)
  (25,0)
  (50,0)
  (75,13.3)
  (100,19.9)
  (125,25.5)
  (150,28.6)
  (175,31.8)
  (200,34.5)
    };

\addplot[
    color=black,
    mark=diamond,
    ]
    coordinates {
  (0,0)
  (25,0)
  (50,23)
  (75,27.4)
  (100,27.4)
  (125,31.3)
  (150,34.6)
  (175,34.8)
  (200,34.8)
    };
    
    \addlegendentry{-} 
    \addlegendentry{0 €/kg} 
    \addlegendentry{2 €/kg}
    \addlegendentry{4 €/kg}
    \addlegendentry{6 €/kg}    
    
    \end{axis}
        \begin{axis}[
        xmin=0,xmax=1,
        xtick={0,0.25,0.5,0.75,1},
        xticklabels={0,0.25M,0.5M,0.75M,1M},
        hide y axis,
        axis x line*=none,
        x tick label style={font=\large},
        ymin=0, ymax=32000,
        x label style={yshift=-0.5cm},
        x tick label style={yshift=-0.5cm}
        ]
    \end{axis} 
\end{tikzpicture} 
 \caption{ Budget Allocation}\label{fig:HydrogenPriceAllocation}
\end{subfigure}
\caption{Effect of hydrogen selling price}\label{fig:HydrogenPriceFigures}
\end{figure}

In Figure~\ref{fig:HydrogenPriceFigures}, the orange line represents the case where the budget is restricted to renewables, excluding the storage option. In that case, operational cost decreases until a certain point (a budget of $0.5$ M€) as the system is supported by increasing renewables, and thus total conventional generation cost and emission penalty decrease. After that point, we observe that the operational cost remains constant, and additional renewables are not integrated into the system. Due to the potential increase in renewable power curtailments arising from the limited capacity of transmission lines, further investments in renewables would increase operational costs. Consideration of hydrogen storage even with no hydrogen market availability (blue line) changes the cost dynamics. After a budget of $0.5$ M€, the operational cost can be decreased as the integration of hydrogen storage saves curtailment costs and adds profit from arbitrage. With the existence of a hydrogen market, when hydrogen can be sold externally (red, green, and black lines), even higher gains are possible, leading to lower operational costs. This underlines the importance of a functioning hydrogen market on cost dynamics. Figure~\ref{fig:HydrogenPriceAllocation} shows that a high hydrogen selling price increases the percentage of the budget allocated for storage.

\subsubsection{Effect of Curtailment Cost}

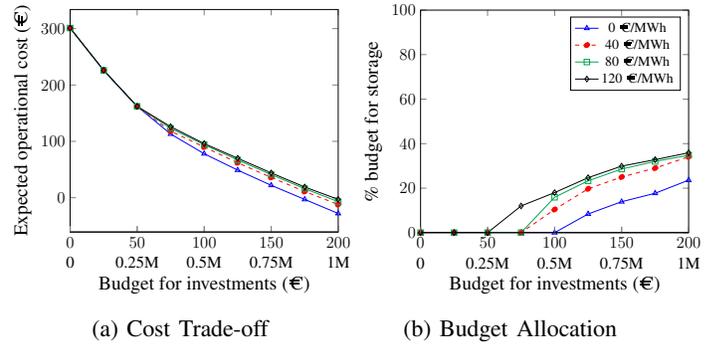
\begin{figure}[h]
\begin{subfigure}{0.25\textwidth}
  \centering
  \begin{tikzpicture}[scale=0.52]
\begin{axis}[
    xlabel={ \Large Budget for investments (€)},
    ylabel={\Large Expected operational cost (€)},
    y tick label style={font=\large},
    x tick label style={font=\large},
    xmin=0, xmax=200,
    xlabel shift={14pt},
    xtick={0,50,100,150,200},
    legend pos=south west,
    grid style=dashed,
]

\addplot[
    color=blue,
    mark=triangle,
    ]
    coordinates {
  (0,301)
  (25,226)
  (50,162)
  (75,113)
  (100,78)
  (125,49)
  (150,22)
  (175,-3)
  (200,-28)
    };

\addplot[
    color=red,
    mark=*,
    style=dashed,
    ]
    coordinates {
  (0,301)
  (25,226)
  (50,162)
  (75,118)
  (100,90)
  (125,62)
  (150,36)
  (175,11)
  (200,-12)
    };

\addplot[
    color=Green,
    mark=square,
    ]
    coordinates {
  (0,301)
  (25,226)
  (50,162)
  (75,123)
  (100,94)
  (125,67)
  (150,41)
  (175,16)
  (200,-7)
    };    
    
\addplot[
    color=black,
    mark=diamond,
    ]
    coordinates {
  (0,301)
  (25,226)
  (50,162)
  (75,126)
  (100,96)
  (125,70)
  (150,44)
  (175,19)
  (200,-3)
    };

\end{axis}
    \begin{axis}[
        xmin=0,xmax=1,
        xtick={0,0.25,0.5,0.75,1},
        xticklabels={0,0.25M,0.5M,0.75M,1M},
        hide y axis,
        axis x line*=none,
        x tick label style={font=\large},
        ymin=0, ymax=32000,
        x label style={yshift=-0.5cm},
        x tick label style={yshift=-0.5cm}
        ]
    \end{axis} 
\end{tikzpicture}

\caption{Cost Trade-off} \label{fig:CurtailmenTrade-off} 
\end{subfigure}
\begin{subfigure}{0.22\textwidth}
  \centering
  \begin{tikzpicture}[scale=0.52]
\begin{axis}[
    xlabel={ \Large Budget for investments (€)},
    ylabel={\Large \% budget for storage},
    y tick label style={font=\large},
    x tick label style={font=\large},
    xmin=0, xmax=200,
    xlabel shift={14pt},
    xtick={0,50,100,150,200},
    ymin=0, ymax=100,
    legend pos=north east,
    grid style=dashed,
]

\addplot[
    color=blue,
    mark=triangle,
    ]
    coordinates {
  (0,0)
  (25,0)
  (50,0)
  (75,0)
  (100,0)
  (125,8.4)
  (150,13.9)
  (175,17.7)
  (200,23.6)
    };

\addplot[
    color=red,
    mark=*,
    style=dashed,
    ]
    coordinates {
  (0,0)
  (25,0)
  (50,0)
  (75,0)
  (100,10.4)
  (125,19.7)
  (150,25)
  (175,29)
  (200,34.1)
  
    };

\addplot[
    color=Green,
    mark=square,
    ]
    coordinates {
  (0,0)
  (25,0)
  (50,0)
  (75,0)
  (100,15.9)
  (125,23.3)
  (150,28.6)
  (175,32)
  (200,34.8)
    };    
    
\addplot[
    color=black,
    mark=diamond,
    ]
    coordinates {
  (0,0)
  (25,0)
  (50,0)
  (75,12)
  (100,18)
  (125,24.7)
  (150,29.9)
  (175,32.8)
  (200,35.9)
    };    
     
    \addlegendentry{0 €/MWh } 
    \addlegendentry{40 €/MWh }
    \addlegendentry{80 €/MWh }
    \addlegendentry{120 €/MWh }

\end{axis}
    \begin{axis}[
        xmin=0,xmax=1,
        xtick={0,0.25,0.5,0.75,1},
        xticklabels={0,0.25M,0.5M,0.75M,1M},
        hide y axis,
        axis x line*=none,
        x tick label style={font=\large},
        ymin=0, ymax=32000,
        x label style={yshift=-0.5cm},
        x tick label style={yshift=-0.5cm}
        ]
    \end{axis} 
\end{tikzpicture}
\caption{Budget Allocation } \label{fig:CurtailmentAllocation} 
\end{subfigure}

\caption{Effect of renewable curtailment cost}\label{fig:CurtailmentFigures}
\end{figure}

Figure~\ref{fig:CurtailmentFigures} depicts that for low investment budgets, changes in curtailment cost do not affect the operational costs since only a few renewables are installed, and hence there is no curtailment. For higher budgets, we can observe this effect since curtailment need arises with the higher penetration of renewables. To prevent a substantial increase in total curtailment costs, the percentage of budget allocated to storage increases, as seen in Figure~\ref{fig:CurtailmentAllocation}. We note that the increase in operational costs would be much more prominent when no storage is available. To examine this further, we solve the model under the same parameter settings with no hydrogen storage. We observe that the network can achieve an average of $34 \%$  operational cost savings with hydrogen storage. This percentage drops to $14 \%$, if we exclude the hydrogen market.

 \subsubsection{Effect of Emission Tax Price}
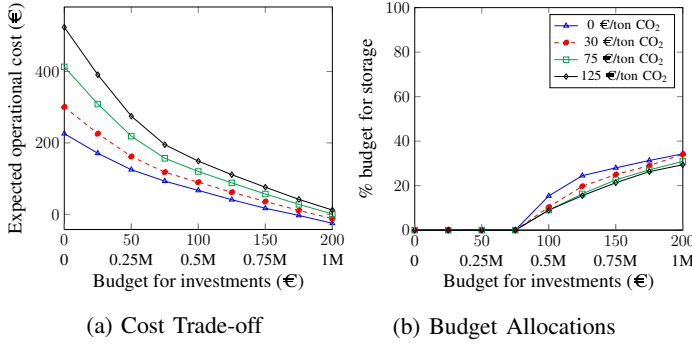
\begin{figure}[h]
\begin{subfigure}{.25\textwidth}
  \centering
  \begin{tikzpicture}[scale=0.52, transform shape]
\begin{axis}[
    xlabel={\Large Budget for investments (€)},
     y tick label style={font=\large},
    x tick label style={font=\large},
    ylabel={\Large Expected operational cost (€)},
    xmin=0, xmax=200,
    xlabel shift={14pt},
    xtick={0,50,100,150,200},
    ymin=-42.5, 
    legend pos=north east,
    ymajorgrids=false,
    grid style=dashed,
]

\addplot[
    color=blue,
    mark=triangle,
    ]
    coordinates {
  (0,226)
  (25,171)
  (50,125)
  (75,93)
  (100,67)
  (125,41)
  (150,17)
  (175,-3)
  (200,-25)
    };
    
\addplot[
    color=red,
    mark=*,
    style=dashed,
    ]
    coordinates {
    (0,301)
  (25,226)
  (50,162)
  (75,118)
  (100,90)
  (125,62)
  (150,36)
  (175,11)
  (200,-12)
    };
    
\addplot[
    color=Green,
    mark=square,
    ]
    coordinates {
  (0,413)
  (25,309)
  (50,219)
  (75,157)
  (100,120)
  (125,88)
  (150,57)
  (175,28)
  (200,1)
    };

\addplot[
    color=black,
    mark=diamond,
    ]
    coordinates {
  (0,524)
  (25,391)
  (50,275)
  (75,195)
  (100,149)
  (125,111)
  (150,76)
  (175,42)
  (200,12)
    };      
\end{axis}
    \begin{axis}[
        xmin=0,xmax=1,
        xtick={0,0.25,0.5,0.75,1},
        xticklabels={0,0.25M,0.5M,0.75M,1M},
        hide y axis,
        axis x line*=none,
        x tick label style={font=\large},
        ymin=0, ymax=32000,
        x label style={yshift=-0.5cm},
        x tick label style={yshift=-0.5cm}
        ]
    \end{axis} 
\end{tikzpicture} 
 \caption{ Cost Trade-off}\label{fig:EmissionTrade-off}
\end{subfigure}
\begin{subfigure}{.22\textwidth}
  \centering
  \begin{tikzpicture}[scale=0.52, transform shape]
\begin{axis}[
    xlabel={\Large Budget for investments (€)},
     y tick label style={font=\large},
    x tick label style={font=\large},
    ylabel={\Large \% budget for storage},
    xmin=0, xmax=200,
    xlabel shift={14pt},
    xtick={0,50,100,150,200},
    ymin=0, ymax=100,
    legend pos=north east,
    ymajorgrids=false,
    grid style=dashed,
]

\addplot[
    color=blue,
    mark=triangle,
    ]
    coordinates {
  (0,0)
  (25,0)
  (50,0)
  (75,0)
  (100,15.4)
  (125,24.5)
  (150,28)
  (175,31.3)
  (200,34.3)
    };
    
\addplot[
    color=red,
    mark=*,
    style=dashed,
    ]
    coordinates {
  (0,0)
  (25,0)
  (50,0)
  (75,0)
  (100,10.4)
  (125,19.7)
  (150,25)
  (175,29)
  (200,34.1)
  
    };
    
\addplot[
    color=Green,
    mark=square,
    ]
    coordinates {
  (0,0)
  (25,0)
  (50,0)
  (75,0)
  (100,9)
  (125,16.3)
  (150,22.6)
  (175,27.1)
  (200,30.9)
    };
    
\addplot[
    color=black,
    mark=diamond,
    ]
    coordinates {
  (0,0)
  (25,0)
  (50,0)
  (75,0)
  (100,9)
  (125,15.5)
  (150,21.3)
  (175,26.3)
  (200,29.4)
    };
    
    \addlegendentry{0 €/ton \ce{\ce{CO2} }} 
    \addlegendentry{30 €/ton \ce{\ce{CO2} }}
    \addlegendentry{75 €/ton \ce{\ce{CO2} }}
    \addlegendentry{125 €/ton \ce{\ce{CO2} }}
                
\end{axis}
    \begin{axis}[
        xmin=0,xmax=1,
        xtick={0,0.25,0.5,0.75,1},
        xticklabels={0,0.25M,0.5M,0.75M,1M},
        hide y axis,
        axis x line*=none,
        x tick label style={font=\large},
        ymin=0, ymax=32000,
        x label style={yshift=-0.5cm},
        x tick label style={yshift=-0.5cm}
        ]
    \end{axis} 
\end{tikzpicture} 
 \caption{Budget Allocations}\label{fig:EmissionAllocation}
\end{subfigure}
\caption{Effect of emission tax price}\label{fig:EmissionFigures}
\end{figure}

Figure~\ref{fig:EmissionTrade-off} shows that the reduction rate of operational cost is diminishing as the investment budget rises under all levels of emission tax price. To examine the effect of including storage, we solve the model under the same parameter settings but with no storage. We observe that the network can achieve an average of $31 \%$  operational cost savings with hydrogen storage. This percentage drops to $16 \%$, if we exclude the hydrogen market. 

Figure~\ref{fig:EmissionAllocation} shows that increasing emission tax may decrease the budget allocated to storage in some specific cases due to the interactions with the hydrogen market. The system owes an emission tax price per kWh of electricity produced by conventional generators. If the renewable output is used at once or stored to meet network demand of a later period, the conventional generation amount decreases. However, if the stored hydrogen is sold outside the network, the conventional generation amount within the network is not altered. However, we should note that replacing green hydrogen with other alternative resources, such as natural gas, significantly reduces emissions. Therefore, to correctly assess emission reductions and to promote green hydrogen production, incentives such as tax credit per emission abated are planned to be given \cite{Abatement} for hydrogen. Thus, considering such policies is likely to change budget allocation dynamics and increase investment in hydrogen storage.

\subsubsection{Effect of Conversion Efficiencies}
\begin{figure}[h]
\begin{subfigure}{.25\textwidth}
  \begin{tikzpicture}[scale=0.52, transform shape]
\begin{axis}[
    xlabel={\Large Budget for investments (€)},
    y tick label style={font=\large},
    x tick label style={font=\large},
    ylabel={\Large Expected operational cost (€)},
    xmin=0, xmax=200,
    xlabel shift={14pt},
    xtick={0,50,100,150,200},
    legend pos=north east,
    ymajorgrids=false,
    grid style=dashed,
]

\addplot[
    color=orange,
    mark=x,
    ]
    coordinates {
  (0,301)
  (25,226)
  (50,162)
  (75,118)
  (100,100)
  (125,96)
  (150,96)
  (175,96)
  (200,96)
    };
    
\addplot[
    color=red,
    mark=*,
    style=dashed,
    ]
    coordinates {
  (0,301)
  (25,226)
  (50,162)
  (75,118)
  (100,90)
  (125,62)
  (150,36)
  (175,11)
  (200,-12)
    };
    
\addplot[
    color=blue,
    mark=triangle,
    ]
    coordinates {
  (0,301)
  (25,226)
  (50,162)
  (75,118)
  (100,88)
  (125,58)
  (150,30)
  (175,2)
  (200,-23)
    };
    
\addplot[
    color=Green,
    mark=square,
    ]
    coordinates {
  (0,301)
  (25,226)
  (50,162)
  (75,118)
  (100,85)
  (125,52)
  (150,21)
  (175,-7)
  (200,-36)
    };

\addplot[
    color=black,
    mark=diamond,
    ]
    coordinates {
  (0,301)
  (25,226)
  (50,162)
  (75,118)
  (100,78)
  (125,42)
  (150,8)
  (175,-23)
  (200,-49)
    };

\end{axis}
    \begin{axis}[
        xmin=0,xmax=1,
        xtick={0,0.25,0.5,0.75,1},
        xticklabels={0,0.25M,0.5M,0.75M,1M},
        hide y axis,
        axis x line*=none,
        x tick label style={font=\large},
        ymin=0, ymax=32000,
        x label style={yshift=-0.5cm},
        x tick label style={yshift=-0.5cm}
        ]
    \end{axis} 
\end{tikzpicture} 
  \caption{Cost trade-off }\label{fig:EfficiencyTrade-offs}
\end{subfigure}
\begin{subfigure}{.22\textwidth}
  \begin{tikzpicture}[scale=0.52, transform shape]
\begin{axis}[
    xlabel={\Large Budget for investments (€)},
     y tick label style={font=\large},
    x tick label style={font=\large},
    ylabel={\Large \% budget for storage},
    xmin=0, xmax=200,
    xlabel shift={14pt},
    xtick={0,50,100,150,200},
    ymin=0, ymax=100,
    legend pos=north east,
    ymajorgrids=false,
    grid style=dashed,
]

\addplot[
    color=orange,
    mark=x,
    ]
    coordinates {
   (0,0)
  (25,0)
  (50,0)
  (75,0)
  (100,0)
  (125,0)
  (150,0)
  (175,0)
  (200,0)
    };
    
\addplot[
    color=red,
    mark=*,
    style=dashed
    ]
    coordinates {
  (0,0)
  (25,0)
  (50,0)
  (75,0)
  (100,10.4)
  (125,19.7)
  (150,25)
  (175,29)
  (200,34.1)
    };
    
\addplot[
    color=blue,
    mark=triangle,
    ]
    coordinates {
  (0,0)
  (25,0)
  (50,0)
  (75,0)
  (100,11.7)
  (125,20.5)
  (150,26.2)
  (175,29.2)
  (200,33.3)
    };
    
\addplot[
    color=Green,
    mark=square,
    ]
    coordinates {
  (0,0)
  (25,0)
  (50,0)
  (75,0)
  (100,14.5)
  (125,22.3)
  (150,26.4)
  (175,30.3)
  (200,33.6)
    };

\addplot[
    color=black,
    mark=diamond,
    ]
    coordinates {
  (0,0)
  (25,0)
  (50,0)
  (75,0)
  (100,18)
  (125,23.7)
  (150,28)
  (175,31)
  (200,34.3)
    };
     \addlegendentry{-}  
    \addlegendentry{0.7-0.5} 
    \addlegendentry{0.8-0.65}
    \addlegendentry{0.9-0.8}
    \addlegendentry{1-1}   
\end{axis}
    \begin{axis}[
        xmin=0,xmax=1,
        xtick={0,0.25,0.5,0.75,1},
        xticklabels={0,0.25M,0.5M,0.75M,1M},
        hide y axis,
        axis x line*=none,
        x tick label style={font=\large},
        ymin=0, ymax=32000,
        x label style={yshift=-0.5cm},
        x tick label style={yshift=-0.5cm}
        ]
    \end{axis} 
\end{tikzpicture} 
  \caption{Budget allocation }\label{fig:EfficiencyAllocation}
\end{subfigure}
    \caption{Effect of conversion efficiencies}
    \label{fig:EfficiencyFigures}
\end{figure}
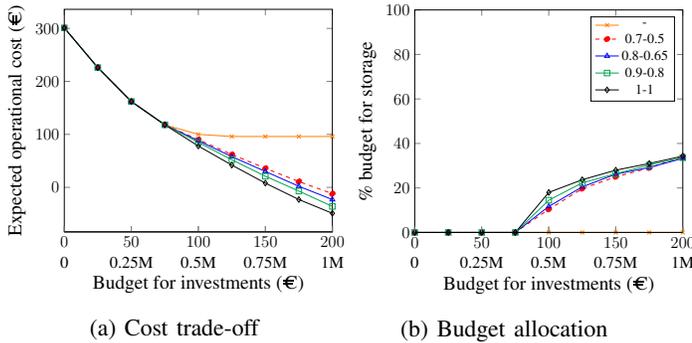

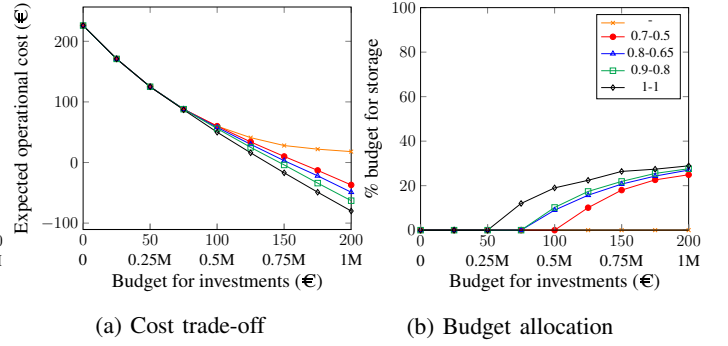
\begin{figure}[h]
\begin{subfigure}{.25\textwidth}
  \begin{tikzpicture}[scale=0.52, transform shape]
\begin{axis}[
    xlabel={\Large Budget for investments (€)},
    y tick label style={font=\large},
    x tick label style={font=\large},
    ylabel={\Large Expected operational cost (€)},
    xmin=0, xmax=200,
    xlabel shift={14pt},
    xtick={0,50,100,150,200},
    legend pos=north east,
    ymajorgrids=false,
    grid style=dashed,
]

\addplot[
    color=orange,
    mark=x,
    ]
    coordinates {
  (0,226)
  (25,171)
  (50,125)
  (75,88)
  (100,60)
  (125,41)
  (150,28)
  (175,22)
  (200,18)
    };
    
\addplot[
    color=red,
    mark=*,
    ]
    coordinates {
  (0,226)
  (25,171)
  (50,125)
  (75,88)
  (100,60)
  (125,34)
  (150,10)
  (175,-13)
  (200,-37)
    };
    
\addplot[
    color=blue,
    mark=triangle,
    ]
    coordinates {
  (0,226)
  (25,171)
  (50,125)
  (75,88)
  (100,58)
  (125,30)
  (150,3)
  (175,-22)
  (200,-49)
    };
    
\addplot[
    color=Green,
    mark=square,
    ]
    coordinates {
  (0,226)
  (25,171)
  (50,125)
  (75,88)
  (100,56)
  (125,25)
  (150,-4)
  (175,-34)
  (200,-63)
    };

\addplot[
    color=black,
    mark=diamond,
    ]
    coordinates {
  (0,226)
  (25,171)
  (50,125)
  (75,87)
  (100,50)
  (125,16)
  (150,-17)
  (175,-49)
  (200,-80)
    };

\end{axis}
    \begin{axis}[
        xmin=0,xmax=1,
        xtick={0,0.25,0.5,0.75,1},
        xticklabels={0,0.25M,0.5M,0.75M,1M},
        hide y axis,
        axis x line*=none,
        x tick label style={font=\large},
        ymin=0, ymax=32000,
        x label style={yshift=-0.5cm},
        x tick label style={yshift=-0.5cm}
        ]
    \end{axis} 
\end{tikzpicture} 
  \caption{ Cost trade-off}\label{fig:EfficiencyTrade-offs2}
\end{subfigure}
\begin{subfigure}{.22\textwidth}
  \begin{tikzpicture}[scale=0.52, transform shape]
\begin{axis}[
    xlabel={\Large Budget for investments (€)},
     y tick label style={font=\large},
    x tick label style={font=\large},
    ylabel={\Large \% budget for storage},
    xmin=0, xmax=200,
    xlabel shift={14pt},
    xtick={0,50,100,150,200},
    ymin=0, ymax=100,
    legend pos=north east,
    ymajorgrids=false,
    grid style=dashed,
]

\addplot[
    color=orange,
    mark=x,
    ]
    coordinates {
   (0,0)
  (25,0)
  (50,0)
  (75,0)
  (100,0)
  (125,0)
  (150,0)
  (175,0)
  (200,0)
    };
    
\addplot[
    color=red,
    mark=*,
    ]
    coordinates {
   (0,0)
  (25,0)
  (50,0)
  (75,0)
  (100,0)
  (125,10.1)
  (150,18)
  (175,22.6)
  (200,24.9)
    };
    
\addplot[
    color=blue,
    mark=triangle,
    ]
    coordinates {
  (0,0)
  (25,0)
  (50,0)
  (75,0)
  (100,9)
  (125,15.7)
  (150,20.7)
  (175,24.3)
  (200,27)
    };
    
\addplot[
    color=Green,
    mark=square,
    ]
    coordinates {
  (0,0)
  (25,0)
  (50,0)
  (75,0)
  (100,10.2)
  (125,17.4)
  (150,21.9)
  (175,25.5)
  (200,27.6)
    };

\addplot[
    color=black,
    mark=diamond,
    ]
    coordinates {
  (0,0)
  (25,0)
  (50,0)
  (75,12)
  (100,19)
  (125,22.4)
  (150,26.4)
  (175,27.4)
  (200,28.9)
    };
    
    \addlegendentry{-}
    \addlegendentry{0.7-0.5} 
    \addlegendentry{0.8-0.65}
    \addlegendentry{0.9-0.8}
    \addlegendentry{1-1}   
\end{axis}
    \begin{axis}[
        xmin=0,xmax=1,
        xtick={0,0.25,0.5,0.75,1},
        xticklabels={0,0.25M,0.5M,0.75M,1M},
        hide y axis,
        axis x line*=none,
        x tick label style={font=\large},
        ymin=0, ymax=32000,
        x label style={yshift=-0.5cm},
        x tick label style={yshift=-0.5cm}
        ]
    \end{axis} 
\end{tikzpicture} 
  \caption{Budget allocation }\label{fig:EfficiencyAllocation2}
\end{subfigure}
    \caption{Effect of conversion efficiencies (curtailment costs and emission tax prices are excluded)}
    \label{fig:EfficiencyFigures2}
\end{figure}

To analyze the effect of conversion efficiencies on decision dynamics, we plot Figures~\ref{fig:EfficiencyFigures}. In Figure~\ref{fig:EfficiencyTrade-offs}, the orange line represents the case without storage. In that case, the operational cost is significantly higher than the cases with storage. This shows that the use of hydrogen storage is cost-efficient in the long term despite the current low conversion efficiencies. As efficiencies improve, we observe that operational costs decreases due to the rise in profit from selling hydrogen and arbitrage revenues. We observe that budget allocation dynamics are not much affected in this particular setting (see Figure~\ref{fig:EfficiencyAllocation}). 

To further elaborate on conversion efficiencies, we exclude curtailment costs and emission tax prices (see Figure~\ref{fig:EfficiencyFigures2}). Compared to Figure~\ref{fig:EfficiencyTrade-offs}, the gap between storage and no-storage options is less in Figure~\ref{fig:EfficiencyTrade-offs2}. In Figure~\ref{fig:EfficiencyAllocation2}, we observe more prominent changes in budget allocation dynamics in comparison to Figure~\ref{fig:EfficiencyAllocation}. We conclude that the curtailment cost and emission tax price can dominate the effect of conversion efficiencies on budget allocation dynamics.

Overall, we observe that daily operational costs are notably higher when we only allow the integration of renewables. It shows the importance of joint optimization of renewables and hydrogen storage integration to achieve operational cost savings. Our optimization framework can provide insights for an investment plan on the economic viability and which part of the investment is made on hydrogen storage to achieve minimum operational cost. Furthermore, our findings show how changes in hydrogen selling price, curtailment cost, emission tax price, and conversion efficiencies affect the budget allocation dynamics. Thus, they can provide valuable insight to authorities on incentivizing network operators to invest in hydrogen storage with regulations in operational level parameters.

\subsection{Optimal Locations and Power Ratings}
In this section, we report corresponding location and power rating decisions. We mainly focus on the cases where we observe significant differences in the budget allocation dynamics in the previous section.

\subsubsection{Effect of Hydrogen Market}
\begin{figure}[h]
\begin{subfigure}{.25\textwidth}
  \begin{tikzpicture}[scale=0.52]
	\begin{axis}[ybar stacked,    
	xlabel={\Large Budget for investments (M€)},
    ylabel={\Large Power rating (kW) },
    ymin=0, ymax=650,
    y tick label style={font=\large},
    x tick label style={font=\large},
    xtick={0,50,100,150,200},
    xticklabels={0,0.25,0.5,0.75,1},
    width=8cm,
    height=6cm,
    legend style={at={(1.9,0.8)},font=\Large}]
	\addplot[color = Green, pattern color=Green, pattern=grid] coordinates
		{(25,109) (50,117) (75,126) (100,135) (125,200) (150,200) (175,200) (200,200) };
	\addplot[color = Goldenrod, pattern color=Goldenrod, pattern=crosshatch] coordinates
		{(25,0) (50,0) (75,0) (100,100) (125,0) (150,0) (175,0) (200,0) };
	\addplot[color = red, pattern color=red, pattern=horizontal lines] coordinates
		{(25,0) (50,100) (75,0) (100,0) (125,0) (150,0) (175,0) (200,0) };
	\addplot[color = orange, pattern color=orange, pattern=north west lines] coordinates
		{(25,0) (50,0) (75,100) (100,0) (125,100) (150,100) (175,100) (200,100) };
	\addplot[color = black, pattern color=black, pattern=dots] coordinates
		{(25,0) (50,0) (75,0) (100,100) (125,100) (150,100) (175,100) (200,100) };
	\addplot[color = blue, pattern color=blue, pattern=vertical lines]  coordinates
		{(25,0) (50,0) (75,100) (100,100) (125,100) (150,100) (175,100) (200,100) };
	\legend{Bus 5,Bus 7,Bus 8,Bus 10,Bus 12,Bus 24}
	\end{axis}
\end{tikzpicture}
  \caption{Renewables (I)}\label{fig:Renewable1}
\end{subfigure}
\newline
\begin{subfigure}{.24\textwidth}
  \begin{tikzpicture}[scale=0.52]
	\begin{axis}[ybar stacked,    
	xlabel={\Large Budget for investments (M€)},
    ylabel={\Large Power rating (kW) },
    ymin=0, ymax=650,
    y tick label style={font=\large},
    x tick label style={font=\large},
    xtick={0,50,100,150,200},
    xticklabels={0,0.25,0.5,0.75,1},
    legend pos=north east,
    width=8cm,
    height=6cm]
	\addplot[color = Green, pattern color=Green, pattern=grid] coordinates
		{(25,109) (50,117) (75,126) (100,181) (125,155) (150,205) (175,221) (200,173) };
	\addplot[color = Goldenrod, pattern color=Goldenrod, pattern=crosshatch] coordinates
		{(25,0) (50,0) (75,0) (100,0) (125,0) (150,0) (175,0) (200,100) };
	\addplot[color = red, pattern color=red, pattern=horizontal lines] coordinates
		{(25,0) (50,117) (75,0) (100,0) (125,100) (150,100) (175,0) (200,0) };
	\addplot[color = orange, pattern color=orange, pattern=north west lines] coordinates
		{(25,0) (50,0) (75,100) (100,115) (125,0) (150,0) (175,119) (200,100) };
	\addplot[color = black, pattern color=black, pattern=dots] coordinates
		{(25,0) (50,0) (75,0) (100,0) (125,100) (150,100) (175,100) (200,100) };
	\addplot[color = blue, pattern color=blue, pattern=vertical lines]  coordinates
		{(25,0) (50,0) (75,100) (100,100) (125,100) (150,100) (175,100) (200,100) };

	\end{axis}
\end{tikzpicture}
  \caption{Renewables (II)}\label{fig:Renewable2}
\end{subfigure}
\begin{subfigure}{.22\textwidth}
  \begin{tikzpicture}[scale=0.52]
	\begin{axis}[ybar stacked,    
	xlabel={\Large Budget for investments (M€)},
    ylabel={\Large Power rating (kW) },
    ymin=0, 
    y tick label style={font=\large},
    x tick label style={font=\large},
    xtick={0,50,100,150,200},
    xticklabels={0,0.25,0.5,0.75,1},
    width=8cm,
    height=6cm]
	\addplot[color = Green, pattern color=Green, pattern=grid] coordinates
		{(25,0) (50,0) (75,0) (100,0) (125,38.1) (150,66.6) (175,76.5) (200,87) };
	\addplot[color = Goldenrod, pattern color=Goldenrod, pattern=crosshatch] coordinates
		{(25,0) (50,0) (75,0) (100,0) (125,0) (150,0) (175,0) (200,30) };
	\addplot[color = red, pattern color=red, pattern=horizontal lines] coordinates
		{(25,0) (50,0) (75,0) (100,0) (125,30) (150,46.8) (175,0) (200,0) };
	\addplot[color = orange, pattern color=orange, pattern=north west lines] coordinates
		{(25,0) (50,0) (75,0) (100,30) (125,0) (150,0) (175,63) (200,49.8) };
	\addplot[color = black, pattern color=black, pattern=dots] coordinates
		{(25,0) (50,0) (75,0) (100,0) (125,0) (150,0) (175,0) (200,30) };
	\addplot[color = blue, pattern color=blue, pattern=vertical lines] coordinates
		{(25,0) (50,0) (75,0) (100,0) (125,0) (150,0) (175,30) (200,30) };
	\end{axis}
\end{tikzpicture}
  \caption{Storage (II)}\label{fig:Storage2}
\end{subfigure}
\newline
\begin{subfigure}{.24\textwidth}
  \begin{tikzpicture}[scale=0.52]
	\begin{axis}[ybar stacked,    
	xlabel={\Large Budget for investments (M€)},
    ylabel={\Large Power rating (kW) },
    ymin=0, ymax=650,
    y tick label style={font=\large},
    x tick label style={font=\large},
    xtick={0,50,100,150,200},
    xticklabels={0,0.25,0.5,0.75,1},
    legend pos=north east,
    width=8cm,
    height=6cm]
    \addplot[color = Green, pattern color=Green, pattern=grid] coordinates
		{(25,109) (50,117) (75,120) (100,169) (125,174) (150,167) (175,207) (200,169) };
	\addplot[color = Goldenrod, pattern color=Goldenrod, pattern=crosshatch] coordinates
		{(25,0) (50,0) (75,0) (100,0) (125,0) (150,0) (175,0) (200,100) };
	\addplot[color = red, pattern color=red, pattern=horizontal lines] coordinates
		{(25,0) (50,100) (75,0) (100,0) (125,0) (150,100) (175,0) (200,0) };
	\addplot[color = orange, pattern color=orange, pattern=north west lines] coordinates
		{(25,0) (50,0) (75,167) (100,179) (125,127) (150,100) (175,112) (200,100) };
	\addplot[color = black, pattern color=black, pattern=dots] coordinates
		{(25,0) (50,0) (75,0) (100,0) (125,0) (150,0) (175,100) (200,100) };
	\addplot[color = blue, pattern color=blue, pattern=vertical lines]  coordinates
		{(25,0) (50,0) (75,0) (100,0) (125,103) (150,100) (175,100) (200,100) };
	\end{axis}
\end{tikzpicture}
  \caption{Renewables (III)}\label{fig:Renewable3}
\end{subfigure}
\begin{subfigure}{.22\textwidth}
  \begin{tikzpicture}[scale=0.52]
	\begin{axis}[ybar stacked,    
	xlabel={\Large Budget for investments (M€)},
    ylabel={\Large Power rating (kW) },
    ymin=0, 
    y tick label style={font=\large},
    x tick label style={font=\large},
    xtick={0,50,100,150,200},
    xticklabels={0,0.25,0.5,0.75,1},
    width=8cm,
    height=6cm
    ]
	\addplot[color = Green, pattern color=Green, pattern=grid] coordinates
		{(25,0) (50,0) (75,0) (100,30.3) (125,39.6) (150,48) (175,60.9) (200,66.3) };
	\addplot[color = Goldenrod, pattern color=Goldenrod, pattern=crosshatch] coordinates
		{(25,0) (50,0) (75,0) (100,0) (125,0) (150,0) (175,0) (200,41.7) };
	\addplot[color = red, pattern color=red, pattern=horizontal lines] coordinates
		{(25,0) (50,0) (75,0) (100,0) (125,0) (150,30) (175,0) (200,0) };
	\addplot[color = orange, pattern color=orange, pattern=north west lines] coordinates
		{(25,0) (50,0) (75,30) (100,36.3) (125,36.9) (150,30) (175,45.3) (200,41.1) };
	\addplot[color = black, pattern color=black, pattern=dots] coordinates
		{(25,0) (50,0) (75,0) (100,0) (125,0) (150,0) (175,41.1) (200,41.1) };
	\addplot[color = blue, pattern color=blue, pattern=vertical lines]  coordinates
		{(25,0) (50,0) (75,0) (100,0) (125,30) (150,35.1) (175,38.7) (200,41.1) };
	\end{axis}
\end{tikzpicture}
  \caption{Storage (III) }\label{fig:Storage3}
\end{subfigure}
    \caption{Power Ratings}
    \label{fig:PowerRatings}
\end{figure}
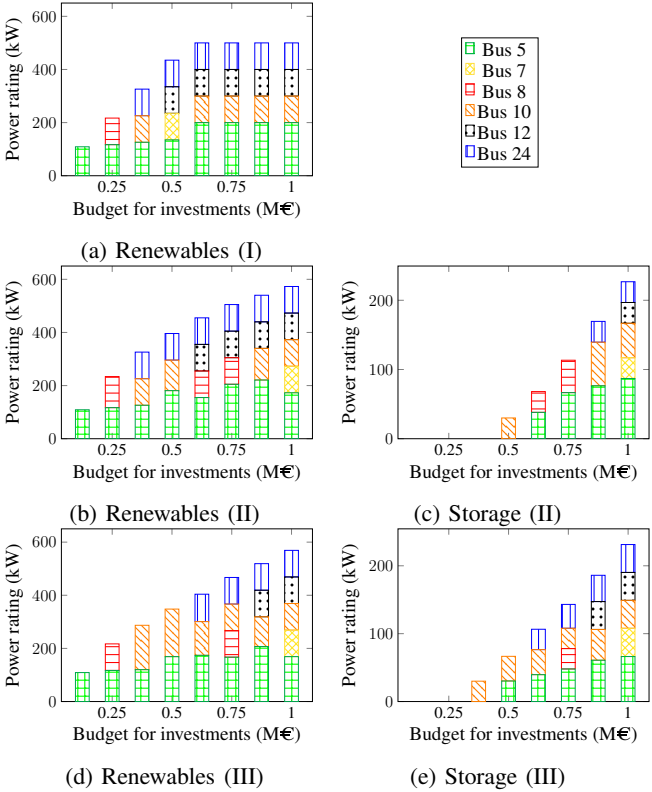

We first analyze the effect of the hydrogen market. Table~\ref{Table:LocationsBudget} shows the change in location decisions for three settings in Figure~\ref{fig:HydrogenPriceFigures}: no-storage, 0 €/kg and 4 €/kg hydrogen selling price. Locations with hydrogen storage are indicated with a superscript plus sign. Figure~\ref{fig:PowerRatings} shows the power ratings of renewables and hydrogen storage for the corresponding locations in Table~\ref{Table:LocationsBudget}.
\begin{table}[ht]
  \centering
    \caption{Locations of Renewables and Hydrogen Storage}
  \label{Table:LocationsBudget}
  \begin{tabular}
  {p{0.08\linewidth}|p{0.15\linewidth}|p{0.28\linewidth}|p{0.31\linewidth}}
       \hline
    Budget & \multicolumn{3}{c}{Locations}  \\
    \hhline{~---}
     (M€) & No-storage (I) & Storage without hydrogen market (II) & Storage with 4 €/kg hydrogen (III)    \\
    \hline
    0.125  & 5          &  5    &  5  \\          
    0.25  & 5 8         &  5 8    &  5 \ 8 \\
    0.375  & 5 10 24    &  5 10 24    &  5 \ 10$^+$ \\
    0.5 & 5 7 12 24     &  5 10$^+$ 24    & 5$^+$ \ 10$^+$ \\
    0.675 & 5 10 12 24  &  5$^+$ 8$^+$ 12 24    & 5$^+$ \ 10$^+$ \ 24$^+$  \\
    0.75& 5 10 12 24    &  5$^+$ 8$^+$ 12 24    & 5$^+$ \ 8$^+$ \ 10$^+$ \ 24$^+$ \\
    0.875 & 5 10 12 24  & 5$^+$ 10$^+$ 12 24$^+$     & 5$^+$ \ 10$^+$ \ 12$^+$ \ 24$^+$ \\
    1 & 5 10 12 24      & 5$^+$ 7$^+$ 10$^+$ 12$^+$ 24$^+$   & 5$^+$ \ 7$^+$ \ 10$^+$ \ 12$^+$  \ 24$^+$ \\
    \hline
  \end{tabular}
\end{table}

For budgets below 0.375 M€, only renewables are located in the same locations in all settings. When the budget is over 0.375 M€, instead of using the entire budget for renewables, we observe a tendency to shift towards hydrogen storage. For example, at a 0.5 M€ budget, allowing for hydrogen storage reallocates the budget from opening renewables at buses 7 and 12 to building hydrogen storage at bus 10 and increasing the power rating of the renewable at bus 5. With a functioning hydrogen market, additional hydrogen storage is built at bus 5, and the renewable at 24 is not opened. When the budget is over 0.675 M€, for the no-storage case, decisions do not change since the increased renewable penetration results in the increased curtailment leading to higher operational costs. For the storage cases, renewables are increasingly co-located with hydrogen storage, increasing total renewable integration in the network. With a hydrogen market, this effect is similar, but the location and sizing decisions in the network differ.

\subsubsection{Effect of the Curtailment Cost}
Table~\ref{Table:Locations2} shows the location decisions corresponding to the 0 and 80 €/MWh curtailment cost settings in Figure~\ref{fig:CurtailmentFigures}.
\begin{table}[h]
  \centering
    \caption{Effect of the Curtailment Cost on Locations}
  \label{Table:Locations2}
  \begin{tabular}
  {p{0.08\linewidth}|p{0.32\linewidth}|p{0.32\linewidth}}
       \hline
    Budget & \multicolumn{2}{c}{Locations}  \\
    \hhline{~--}
     (M€) &  0 €/MWh curtailment cost & 80 €/MWh curtailment cost    \\
    \hline
    0.125  & 5                        & 5     \\          
    0.25   & 5 8                      & 5 10    \\
    0.375  & 5 10                     & 5 10 24     \\
    0.5    & 5 8 21                   & 5 10 24$^+$   \\
    0.675  & 5 8$^+$ 15 27            & 5$^+$ 10$^+$ 24$^+$  \\
    0.75.  & 5$^+$ 7 8 10 15$^+$      & 5$^+$ 8 10$^+$ 24$^+$   \\
    0.875  & 5$^+$ 7 8 15$^+$ 28$^+$  & 5$^+$ 10$^+$ 12$^+$ 24$^+$\\
    1      & 5$^+$ 7 8$^+$ 12 24$^+$  & 5$^+$ 7$^+$ 10$^+$ 12$^+$ 24$^+$  \\
    \hline
  \end{tabular}
\end{table}

Regardless of the change in parameters, bus 5 is always a preferred location that has the highest power load. For 0 €/MWh curtailment costs, we also observe buses 7 and 8, which are the locations with the highest demand after bus 5. For 80 €/MWh curtailment costs, buses 10 and 24 are frequently preferred. Compared to buses 7 and 8, the total thermal limit of the transmission lines connected to bus 10 is higher, which accommodates renewables conveniently by dispatching excess power. We observe that an increase in curtailment cost shifts location from higher power loads to higher thermal limits.

Particularly the locations near load centers are preferable since they reduce supply needs from more distant generators, thereby reducing transmission losses. Transmission lines connected to them have higher total thermal limits, so curtailment need is less. The resulting location decisions show the importance of considering transmission losses and, thereby, the importance of AC power equations.

\subsubsection{Comparison with DC Approximation}
We obtain the location and sizing decisions from the DC approximation of the original model rather than its MISOCP relaxation. Then, we fix the decisions in the original AC formulation to make operational cost comparisons. The results are reported in Table~\ref{Table:AC/DC} for the base case with red dashed lines. We note, on average, 23\% operational cost savings with the AC formulation. The results emphasize the importance of considering the AC OPF dynamics for economically efficient location and sizing decisions. 

\begin{table}[h]
  \centering
    \caption{Comparison with DC Approximation}
  \label{Table:AC/DC}
  \begin{tabular}  {p{0.09\linewidth}|p{0.055\linewidth}|p{0.045\linewidth}|p{0.055\linewidth}|p{0.045\linewidth}|p{0.055\linewidth}|p{0.045\linewidth}|p{0.055\linewidth}|p{0.045\linewidth}|p{0.04\linewidth}}
   \hline 
 Budget &0.125 & 0.25 & 0.375 & 0.5 & 0.675 & 0.75 & 0.875 & 1 & Avg \\
 \hline
  AC & 226 & 162 & 118 & 89 & 62 & 36 & 11 & -12 & 86 \\
    DC & 229 & 167 & 134 & 111 & 89 & 71 & 56 & 38 & 112 \\
    Gap \% & 1  & 3  & 12  & 20  & 30  & 48  & 80 & 131 & 23 \\
    \hline
  \end{tabular}
\end{table}

\section{Conclusions}\label{sec:conclusions}
This paper proposes a joint optimization model for the location and sizing of renewables and hydrogen storage based on multi-period AC OPF. We provide a systematic solution approach based on SOCP within a Benders decomposition framework to provide solutions to our model with a global optimality guarantee. On a representative test case, we conduct computational experiments and show that the joint integration of renewables and hydrogen storage leads to significant operational cost savings. Furthermore, we show that it is crucial to consider AC power flow equations instead of DC approximations as they lead to the improved location and sizing decisions and thus lower operational costs. Moreover, we show how a functioning hydrogen market can change decision dynamics. Finally, we use our solution framework to provide qualitative insights for decision-makers on how to integrate renewables and hydrogen storage under varying operational parameters such as the hydrogen selling price, curtailment cost, emission tax price, and conversion efficiency.

Our optimization framework is general, meaning that the operational specifications of MOPF, investment decisions, and storage type can be adapted or altered without affecting the structural ideas of our solution method. 

Future research might focus on considering uncertainties in future hydrogen markets regarding prices and demands. A natural next step for our research is considering the expansion of the transmission lines in the network too.

\section*{Acknowledgment}
This project has received funding from the Fuel Cells and Hydrogen 2 Joint Undertaking under grant agreement No 875090, HEAVENN - Hydrogen Energy Applications in Valley Environments for Northern Netherlands. 

\ifCLASSOPTIONcaptionsoff
  \newpage
\fi

\bibliographystyle{ieeetr}
\bibliography{references}

\end{document}